\documentclass[12pt,twoside,amstex]{article}
\usepackage{amsmath,amsfonts}
\def\R{{\mathord{\rm I\mkern-3.6mu R}}}
\def\N{{\mathord{\rm I\mkern-3.6mu N}}}
\def\Z{{\mathord{\rm Z\mkern-5.6mu Z}}}
\newcommand{\me}{1/2}
\newcommand{\D}{\protect\displaystyle}
\newcommand{\T}{\protect\textstyle}
\newcommand{\ipl}{\langle}
\newcommand{\ipr}{\rangle}
\newcommand{\Lbd}{\Lambda}
\newcommand{\lbd}{\lambda}
\newcommand{\eps}{\varepsilon}
\newcommand{\vphi}{\varphi}

\begin{document}

\title{\centering \large \bf MEAN VALUE METHODS FOR SOLVING THE HEAT 
EQUATION BACKWARDS IN TIME}
\author{\large \bf A. Leit\~ao\thanks{Partially supported by CNPq grant 
150019/01-7.}}
%
%
\date{}
\maketitle

\vspace{-6.1cm}
\noindent {\texttt{Proc.\,of\;VII\:Workshop\:on\:PDE's\;/\;%
July\:16\,-\,20,\:2001\;/\;Rio\:de\:Janeiro,\:Brazil}
\vspace{5.2cm}

\thispagestyle{empty}
\pagestyle{myheadings}
\markboth{A. Leit\~ao}
  {Mean value methods for solving the heat equation ...}

\setlength{\baselineskip}{6mm}
%
%

\newtheorem{definition}{Definition}
\newtheorem{example}{Example}
\newtheorem{algorithm}{Algorithm}
\newtheorem{application}{Application}
\newtheorem{theorem}{Theorem}
\newtheorem{exercise}{Exercise}
\newtheorem{lemma}{Lemma}
\newtheorem{sublemma}{Sublemma}
\newtheorem{notation}{Notation}
\newtheorem{observation}{Observation}
\newtheorem{comment}{Comment}
\newtheorem{consequence}{Consequence}
\newtheorem{illustration}{Illustration}
\newtheorem{problem}{Problem}
\newtheorem{proposition}{Proposition}
\newtheorem{proposal}{Proposal}
\newtheorem{property}{Property}
\newtheorem{remark}{Remark}
\newtheorem{hint}{Hint}
\newtheorem{result}{Result}
\newtheorem{conjecture}{Conjecture}
\newtheorem{claim}{Claim}
\newtheorem{note}{Note}
\newtheorem{question}{Question}
\newtheorem{fact}{Fact}
\newtheorem{assumption}{Assumption}
\newtheorem{corollary}{Corollary}
\newtheorem{principle}{Principle}
%
%
\def\Box{ \framebox[5,5pt]}

\begin{abstract}
We investigate an iterative mean value method for the inverse (and highly 
ill-posed) problem of solving the heat equation backwards in time. 
Semi-group theory is used to rewrite the solution of the inverse problem 
as the solution of a fixed point equation for an affine operator, with
linear part satisfying special functional analytical properties. 
We give a convergence proof for the method and obtain convergence rates 
for the residual. Convergence rates for the iterates are also obtained 
under the so called source conditions.
\end{abstract}
%
%
%
\section{Introduction} \label{sec1}

The problem of solving the heat equation backwards in time is a classical 
example of ill-posed problem. Using semi-group theory to represent the 
solution of parabolic problems, one can verify that this particular problem 
is modeled by a linear positive operator, whose eigenvalues converge 
exponentially to zero (see Section~\ref{sec2}).

In this paper we introduce a mean value method, based on the Mann iteration, 
in order to find the solution of a fixed point equation associated to the 
inverse problem. The development presented in [6] for elliptic Cauchy 
problems is extended to this inverse parabolic problem.

Our algorithm generalizes the iteration introduced in [1], where this 
problem is also considered. In that paper the authors obtain the same fixed 
point equation treated here. However they propose a method based on solving 
successive well-posed initial value problems, which can be interpreted as a 
particular case of our mean value method.

The paper is organized as follows. In Section~\ref{sec2} we introduce some 
necessary notation and formulate the inverse problem. We also reinterpret 
the inverse problem in the form of a fixed point equation. In Section~%
\ref{sec3} we present an overview of the iterative method introduced by 
W.\,Mann. In Section~\ref{sec4} we formulate our mean value method and 
give a convergence proof for the case of exact data. In Section~\ref{sec5} 
we consider problems with noisy data. The generalized discrepancy principle 
is used to obtain convergence rates for the iteration residual. Using 
appropriate {\em source conditions}, we also prove convergence rates for 
the iterates. Finally, we discuss an example, which shows that for this 
particular inverse problem, the source conditions can be interpreted in 
terms of regularity of the solution in Sobolev spaces.
%
%
%
\section{Formulation of the inverse problem} \label{sec2}

We start this section introducing some notation. Given a normed linear 
space $H$, we call an operator $T: H \to H$ {\em non expansive} if 
$\|T\| \leq 1$. An arbitrary operator $T: H \to H$ is called {\em regular 
asymptotic} in $H$ if
$$ \lim_{k\to\infty} \| T^{k+1}(x) - T^k(x) \| \to 0 ,\ \forall x \in H . $$

Let $\Omega \subset \R^n$ be an open bounded set with smooth boundary 
$\partial\Omega$ and $\Lbd$ be a positive self-adjoint unbounded operator 
(with discrete spectrum) densely defined on $L^2(\Omega)$. The resolution 
of the identity associated to $\Lbd$ is represented by $E_\lbd$, $\lbd 
\in \R$. The family of Hilbert spaces ${\cal H}^s(\Omega)$, $s \geq 0$ (or 
simply ${\cal H}^s$) is now defined by
\begin{equation} \label{def:Hs-raeume}
{\cal H}^s(\Omega) \ := \ \{ \vphi \in L^2(\Omega); \ \| \vphi \|_s^2 := 
\int_0^\infty (1+\lbd^2)^s d \ipl E_\lbd \vphi, \vphi \ipr
  < \infty \} .
\end{equation}
Obviously  ${\cal H}^0 = L^2(\Omega)$. The Hilbert spaces 
${\cal H}^{-s}(\Omega)$ (with $s > 0$) are defined by duality: 
${\cal H}^{-s} := ({\cal H}^{s})'$. One should notice that the 
embedding ${\cal H}^r \hookrightarrow {\cal H}^s$, $r > s$ is 
dense and compact. It is worth mentioning that in the special 
case $\Lbd = (-\Delta)^{\me}$, where $\Delta$ is the Laplace--Beltrami 
operator on $\Omega$, we have ${\cal H}^s(\Omega) = H^{2s}_0(\Omega)$, 
where $H^{s}_0(\Omega)$ are the Sobolev spaces according to Lions 
and Magenes.

Next we formulate the inverse problem considered in this paper. Given a 
function $f \in {\cal H}^0$ let $u \in (V, \|\cdot\|_{V})$, where
\begin{equation}
V := L_2(0,T; {\cal H}^1(\Omega)) ,\ \ \
\|u\|_V^2 := \int_0^T ( \|u(t)\|_1^2 + \|\partial_t u(t)\|_{-1}^2 )\, dt ,
\end{equation}
be the solution of \\[2ex]
$ (P) \hfill
      \left\{  \begin{array}{l}
         (\partial_t + \Lbd^2) u \, = \, 0 ,\ {\rm in}\ (0,T) \times \Omega \\
         u(T)               \, = \, f\, .
      \end{array} \right. \hfill  $ \\[2ex]
Our goal is to reconstruct $u(t)$ at $t=0$. In other words, given a 
temperature profile at the time $t=T$ and a heat transport equation, 
find the corresponding temperature profile at the initial time $t=0$. 
Problem $(P)$ corresponds to the inverse heat transport problem, which 
is known to be exponentially ill-posed. This assertion follows 
immediately from the explicit representation of the solution $u$ of 
problem $(P)$, which is given by
\begin{equation} \label{eq:loesung-parab}
  u(t) \ = \ \exp(\Lbd^2 (T-t)) \, f .
\end{equation}
Indeed, using the eigenfunctions of $\Lbd$, one can construct a sequence 
of data $f_k$ converging (uniformly) to zero in $H$, while the $V$-norm 
of the corresponding solutions $u_k$ do not.

Notice that if $u \in V$, then $u(0),\ u(T) \in {\cal H}^0$. Another 
important remark concerning the inverse problem: given $f \in 
{\cal H}^0$, problem $(P)$ has exactly one solution in $V$. Thus 
the problem of determining $u(0)$ (solution) from $u(T)$ (data) 
has always a solution, this solution is unique, but it does not 
depend in a stable way on the problem data.

In the sequel we characterize the solution of our inverse problem as the 
solution of a fixed point equation. We start by considering problem $(P)$ 
with data $f \in {\cal H}^0$ and denote by $\bar{x}$ the solution of 
the inverse problem, i.e. $\bar{x} = u(0)$. Let the positive constant 
$\bar{\lbd}$ be defined by
\begin{equation} \label{eq:gammaB-wahl}
\bar{\lbd} := \inf\{ \lbd;\ \lbd \in \sigma(\Lbd) \} .
\end{equation}
Next we chose a parameter
\begin{equation} \label{eq:gamma-wahl}
\gamma \in (0 , 2 \exp(\bar{\lbd}^2 T)) .
\end{equation}
Given $\vphi \in {\cal H}^0$, let us consider the following initial value 
problem of parabolic type \\[2ex]
$ (Q) \hfill
      \left\{  \begin{array}{l}
        (\partial_t + \Lbd^2) w = 0,\ {\rm in}\ (0,T) \times \Omega \\
        w(0) = \vphi
      \end{array} \right. \hfill $ \\[2ex]
We define the affine operator $T: {\cal H}^0 \to {\cal H}^0$ by
\begin{equation} \label{eq:T-def1}
T \vphi := \vphi - \gamma (w(T) - f) .
\end{equation}
A straightforward calculation shows that $T \bar{x} = \bar{x}$, 
i.e. $\bar{x}$ is a fixed point of $T$. Further, since the solution 
of problem $(Q)$ can be written as
\begin{equation}
w(t) = \exp(-\Lbd^2 t)\, \vphi ,
\end{equation}
we obtain for the operator $T$ the representation
\begin{equation} \label{eq:T-def2}
T \vphi = \big( I - \gamma \exp(-\Lbd^2 T) \big)\, \vphi + \gamma f .
\end{equation}

\begin{remark} \label{rem:Tl}
In [1] several properties of the linear part of the operator $T$, 
namely $T_l = I - \gamma \exp(-\Lbd^2 T)$, are investigated. The most 
relevant ones are self-adjointness, nonexpansivity, asymptotic 
regularity and the fact that 1 is not an eigenvalue of $T_l$ (although 
1 belongs to the continuous spectrum of $T_l$). Further, under the 
stronger assumption $\gamma < 2 \exp(\tilde{\lbd}^2T)$, where 
$\tilde{\lbd} := ( \bar{\lbd}^2 - T^{-1} \ln 2 )^{\me}$, injectivity 
of $T_l$ can also be proved.
\end{remark}

\begin{remark}
If we write the inverse problem in the form $S \bar{x} = f$, 
with $S := \exp(-\Lbd^2 T)$, equation \eqref{eq:T-def1} resembles very 
much the fixed point equation
\begin{equation}
\bar{x} = \bar{x} + S^* (f - S \bar{x}) ,
\end{equation}
which is based on a transformation of the normal equation. If 
$\|S\|^2 < 2$, the related fixed point operator $I - S^*S$ is 
nonexpansive and one may apply the method of successive approximation 
(this corresponds to the so called Landweber iteration: $\vphi_{k+1} = 
\vphi_k + S^* (f - S\vphi_k)$). However, this condition is not 
satisfied is our case. The idea in [1] is to choose a relaxation 
parameter $\gamma >0$ such that $I - \gamma S$ is nonexpansive 
and then apply the successive approximation method for the fixed point 
equation
\begin{equation}
\bar{x} = \bar{x} + \gamma\, (f - S \bar{x}) .
\end{equation}
\end{remark}
%
%
%
\section{The Mann iteration} \label{sec3}

We present a brief overview of the method introduced by W.\;Mann in 1953 
in [9]. Given a Banach space $X$ and $E \subset X$, Mann considered 
the problem of approximating a solution of the fixed point equation for 
a continuous operator ${\cal T}: E \to E$.

In order to avoid the problem of existence of fixed points, the subset 
$E$ was assumed to be convex and compact (the existence question is than 
promptly answered by the Schauder fixed point theorem). Strongly influenced 
by the works of Ces\`aro and Topelitz, who used mean value methods in the 
summation of divergent series, Mann proposed a mean value iterative method 
based on the Picard iteration ($x_{k+1} := {\cal T}(x_k)$), which we shall present 
next.

Let $A$ be the (infinite) lower triangular matrix
$$ A \, = \, \left( \begin{array}{cccccc}
     1      & 0      & 0      & \cdots & 0      & \cdots \\
     a_{21} & a_{22} & 0      & \cdots & 0      & \cdots \\
     a_{31} & a_{32} & a_{33} & \cdots & 0      & \cdots \\
     \vdots & \vdots & \vdots & \ddots & \vdots & \vdots
   \end{array} \right) , $$
with coefficients $a_{ij}$ satisfying 
\begin{itemize}
\item[\it i)]   $a_{ij} \, \ge \, 0\, ,\ i,j = 1, 2, \ldots$;
\item[\it ii)]  $a_{ij} \, = \, 0\, ,\ j > i$;
\item[\it iii)] $\sum\limits_{j=1}^i a_{ij} \, = \, 1\, , \ 
                i = 1, 2, \ldots$.
\end{itemize}
Starting with an arbitrary element $x_1 \in E$, the {\em Mann iteration} is 
defined in the following way
\begin{tt}
\begin{itemize} \itemsep0.1ex
\item[1.] Choose $x_1 \in E$;
\item[2.] For $k = 1, 2, \ldots$ do \\[0.5ex]
     \mbox{\ \ } $v_k := \sum_{j=1}^k \, a_{kj} \, x_j$; \\[1ex]
     \mbox{\ \ } $x_{k+1} := {\cal T}(v_k)$;
\end{itemize}
\end{tt}
In the sequel we denote this iteration briefly by $M(x_1, A, {\cal T})$. Notice 
that with the particular choice $A = I$, this method corresponds to the 
usual successive approximation method. Next we state the main result 
in [9].

\begin{lemma} \label{lemma:mann}

Let $X$ be a Banach space, $E \subset X$ a convex compact subset, 
${\cal T}: E \to E$ continuous. Further, let $\{ x_k \}$, $\{ v_k \}$ 
be the sequences generated by the iteration $M(x_1, A, {\cal T})$. If 
either of the above sequences converges, then the other also converges 
to the same point, and their common limit is a fixed point of ${\cal T}$.
\end{lemma}

\noindent{\bf Proof.} See Theorem~1 in [9].

\hfill \Box

In that paper, the case where neither of the sequences $\{ x_k \}$, 
$\{ v_k \}$ converges is also considered. Under additional requirements 
on the coefficients $a_{ij}$, a relation between the sets of limit 
points of $\{ x_k \}$ and $\{ v_k \}$ is proven.

Many authors considered the Mann iteration in other frameworks. Among 
others we mention [4], [7], [8] and [10]. In the sequel we discuss 
a result obtained by C.\,Groetsch, which will be useful in the 
further discussion. In [7], Groetsch considers a variant of the Mann 
iteration. The Matrix $A$ is assumed to be {\em segmenting}, i.e. 
additionally to properties \,{\it i)}, \,{\it ii)} and \,{\it iii)}, 
the coefficients $a_{ij}$ must also satisfy
\begin{itemize}
\item[\it iv)] $a_{i+1,j} \, = \, (1 - a_{i+1,i+1}) \, a_{ij}$, \ $j \le i$.
\end{itemize}
One can easily check that, under assumptions {\it i)}, \dots, {\it iv)}, 
the element \,$v_{k+1}$ of $M(x_1, A, {\cal T})$ can be written in the 
form of the convex linear combination
\begin{equation} \label{eq:segment}
v_{k+1} \ = \ (1 - d_k) v_k \, + \, d_k {\cal T}(v_k) ,
\end{equation}
where $d_k := a_{k+1,k+1}$. In other words, $v_{k+1}$ lies on the line 
segment joining $v_k$ and $x_{k+1} = {\cal T}(v_k)$, what justifies the 
denomination of {\em segmenting matrix}. Notice that in this case the 
choice of the diagonal elements $d_k$ determines completely the matrix 
$A$. The following lemma corresponds to the main result in [7].

\begin{lemma} \label{lemma:groe}

Let $X$ be an {\em uniformly convex} Banach space, $E \subset X$ a convex 
subset, ${\cal T}: E \to E$ a nonexpansive operator with at least one fixed 
point in $E$. If $\sum_{k=1}^\infty \, d_k (1-d_k)$ diverges, then the 
sequence $\{ (I-{\cal T}) v_k \}$ converges strongly to zero, for every 
$x_1 \in E$.
\end{lemma}

\noindent{\bf Proof.} See Theorem~2 in [7].

\hfill \Box

In order to prove strong convergence of the sequence $\{ x_k \}$, one needs 
stronger assumptions on both the set $E$ and the operator ${\cal T}$ (e.g. 
$E$ is also closed and ${\cal T}(E)$ is relatively compact in $X$).

Notice that Lemma~\ref{lemma:groe} gives on $\{ x_k \}$ a condition 
analogous to the {\em asymptotic regularity}, which is used in both [10] 
and [2]. This condition is also used in [1] and [6] for the analysis of 
linear Cauchy problems.

The last result we analyze in this section is due to H.W.\,Engl and 
A.\,Leit\~ao and is the analog of Lemma~\ref{lemma:groe} for affine 
operators with nonexpansive linear part, defined on Hilbert spaces.

\begin{lemma} \label{lemma:enle}
Let $H$ be a Hilbert space and ${\cal T}: H \to H$ an affine operator 
with nonexpansive linear part. Further, let $A$ be a segmenting matrix 
such that $\sum_{k=1}^\infty \, d_k (1-d_k)$ diverges. The iteration 
$M(x_1, A, {\cal T})$ generates a sequence $v_k$ such that 
$\{ (I - {\cal T}) v_k \}$ converges strongly to zero, for 
every $x_1 \in H$.
\end{lemma}

\noindent{\bf Proof.} See Theorem~6 in [6].

\hfill \Box

\section{A mean value method for the inverse heat transport problem}
\label{sec4}

The iterative method proposed in this paper corresponds to the Mann 
iteration applied to the fixed point equation $T \vphi = \vphi$, where 
the operator $T$ is defined in~\eqref{eq:T-def2}. 
Initially we address the question of convergence for exact data. Given 
the exact data $f \in {\cal H}^0$ we have the following convergence result:

\begin{theorem}
Let $T$ be the operator defined in~\eqref{eq:T-def2} and $A$ a 
segmenting matrix with $\sum_{k=1}^{\infty} d_k (1-d_k) = \infty$. 
For every $x_1 \in {\cal H}^0$ the iteration $M(x_1, A, T)$ generates 
sequences $x_k$ and $v_k$, which converge strongly to $\bar{x}$, 
the uniquely determined fixed point of $T$.
\end{theorem}

\noindent{\bf Proof.}
Existence and uniqueness of the fixed point $\bar{x}$ was already 
justified on Section~\ref{sec2}. Let $T_l$ be the linear part of the 
operator $T$. Since $(I - T_l) (v_k - \bar{x}) = (I - T) v_k$ and 
Ker$(I - T_l) = \{ 0 \}$ (see Remark~\ref{rem:Tl} or [1, Lemma~8] 
for details), it is enough to prove that $\lim_k (I - T) v_k = 0$. 
This however follows from Lemma~\ref{lemma:enle}.

\hfill \Box

%
%
%
\section{Convergence rates} \label{sec5}

In this section we consider the case of inexact data. This is a 
particularly interesting question, what concerns the application of 
this method to {\em real live problems}. Let us assume we are given 
noisy Cauchy data $f_\eps \in {\cal H}^0$, such that
\begin{equation} \| f - f_\eps \| \leq \eps , \end{equation}
where $\eps>0$ is the noise level and $f \in {\cal H}^0$ represents the 
exact problem data. This assumption is relevant in the case where $f_\eps$ 
is obtained by means of measurements, since measured data always contain 
errors. In this case we consider the iteration residual and use the 
generalized {\em discrepancy principle} to provide a stopping rule for 
the algorithm.

In the sequel, we consider for simplicity the particular case where the 
segmenting matrix is given by $A=I$, the identity matrix. We start by 
defining the iteration residual. Given $x_1 \in {\cal H}^0$, let us 
consider the sequences
\begin{eqnarray}
x_{k+1}      & = & T_l\, x_k + \gamma f , \\
x_{k+1}^\eps & = & T_l\, x_k^\eps + \gamma f_\eps ,
\end{eqnarray}
generated by the iterative method (note that $x_1^\eps = x_1$). The 
corresponding residuals (exact and real, i.e. using noisy data) are 
defined by
\begin{eqnarray}
r_k      & := & \gamma f - (I-T_l) x_k , \\
r_k^\eps & := & \gamma f_\eps - (I-T_l) x_k^\eps .
\end{eqnarray}
Notice that the residual sequences $\{ \| r_k^\eps \|\}$, $\{ \| r_k \|\}$ 
are nonincrea\-sing. Indeed, this follows from the nonexpansivity of the 
operator $T_l$.

Now let $\mu > 1$ be fixed. According to the discrepancy principle, the 
iteration should be stopped at the step $k(\eps,f_\eps)$, when for the 
first time $ \| r_{k(\eps,f_\eps)}^\eps \| \le \mu \eps$, i.e.
\begin{equation} \label{eq:discrep}
k(\eps,f_\eps) \ := \ \min \{ k \in \N\ |\ 
   \| \gamma f_\eps - (I-T_l) \vphi_k^\eps \| \le \mu\eps \} .
\end{equation}
If we do not make any further (regularity) assumption on the solution 
$\bar{x}$, we cannot prove convergence rates for the iterates 
$\| x_k - \bar{x} \|$. However, it is possible to obtain rates 
of convergence for the residuals.

Next we obtain an estimate for $k(\eps,f_\eps)$ in \eqref{eq:discrep}. 
The proof of this result is similar to the one known for the {\em Landweber 
iteration} (see, e.g. [5, Section~6.1]). For convenience of the reader 
we include here the proof.

\begin{proposition} \label{prop:rate-discr}
If $\mu >1$ is fixed, the stopping rule defined by the {\em discrepancy 
principle} in \eqref{eq:discrep} satisfies $k(\eps,f_\eps) = O(\eps^{-2})$.
\end{proposition}

\noindent{\bf Proof.}
We start from the identity
\begin{multline} \label{eq:prop-rd1}
\| \bar{x} - x_{j+1} \|^2 \, = \, \| \bar{x} - x_j \|^2
 - \| \gamma f - (I-T_l) x_j \|^2 \\
 - 2 \ipl (I-T_l) (\bar{x} - x_j), T_l (\bar{x} - x_j) \ipr
\end{multline}
to obtain the estimate
\begin{equation} \label{eq:prop-rd2}
\| \bar{x} - x_j \|^2 - \| \bar{x} - x_{j+1} \|^2
  \ge \| \gamma f - (I-T_l) x_j \|^2 .
\end{equation}
Adding up this inequalities for $j = 1, \ldots, k$ we can conclude
\begin{equation} \label{eq:prop-rd3}
\| \gamma f - (I-T_l) \vphi_k \|^2 \le k^{-1} \| \bar{x} - x_1 \|^2 .
\end{equation}
Since the real residual $r_k^\eps$ satisfy
\begin{equation} \label{eq:prop-rd4}
\| \gamma f_\eps - (I-T_l) x_{k+1}^\eps \|
 \le \eps + \| \gamma f - (I-T_l) x_k \| ,
\end{equation}
we obtain from \eqref{eq:prop-rd3}
\begin{equation} \label{eq:prop-rd5}
\| \gamma f_\eps - (I-T_l) x_{k+1}^\eps \| \le
   \eps \, + \, k^{-\frac{1}{2}} \, \| \bar{x} - x_1 \| .
\end{equation}
Since the right hand side of \eqref{eq:prop-rd5} is lower than $\mu\eps$ 
for $k > (\mu-1)^{-2} \| \bar{x} - x_1 \|^2 \eps^{-2}$, we have 
$k(\eps,z_\eps) \le c\, \eps^{-2}$, where the constant $c > 0$ depends 
only on $\mu$ and $x_1$.

\hfill \Box

From Proposition~\ref{prop:rate-discr} it is possible to obtain rates 
of convergence for the residuals.

\begin{corollary}
Let $f$ be the exact data, $\tau > 1$ and $\eps > 0$. Given
noisy data $f_\eps$, with $\| f_\eps - f \| \le \eps$, the
stopping rule $k(\eps,f_\eps)$ determined by the discrepancy
principle satisfies
\begin{itemize}
\item [i)]  $\| \gamma f_\eps - (I-T_l) x_{k(\eps,f_\eps)}^\eps \| \le 
            \mu \eps$;
\item [ii)] $k(\eps,z_\eps) = O(\eps^{-2})$.
\end{itemize}
\end{corollary}

If appropriate regularity assumptions are made on the fixed point 
$\bar{x}$, it is possible to obtain convergence rates also for 
the approximate solutions. This additional assumptions are stated 
here in the form of the so-called {\em source conditions} (see, e.g., 
[5]). Since our inverse problem is exponentially ill-posed, the 
source condition take the form
\begin{equation} \label{gl:source-cond}
\bar{x} - x_1 = F(I - T_l) \, y ,
\end{equation}
where $y$ is some function in ${\cal H}^0$ and $F$ is defined by
$$  F(\lbd) := \left\{ \begin{array}{cl}
                 (\ln(e/\lbd))^{-p}, & \lbd > 0 \\
                                  0, & \lbd = 0
               \end{array} \right. $$
with $p>0$ fixed. This choice of $F$ corresponds to the {\em 
logarithmic-type source conditions}. Under these assumptions we 
can prove the following rates:

\begin{proposition} \label{prop:source}
Let $f$ be given data and assume that the fixed point $\bar{x}$ of 
$T$ satisfies the source condition
\begin{equation} \label{eq:source-cond}
\bar{x} - x_1 \; = \; F(I-T_l)\, y , \mbox{ for some }
                       y \in {\cal H}^0 ,
\end{equation}
where $x_1 \in {\cal H}^0$ is some initial guess and $F$ is defined 
as above for some $p \ge 1$. Let $\mu > 2$, $f_\eps$ some given 
noisy data with $\| f_\eps - f \| \le \eps$, $\eps > 0$ and 
$k(\eps,f_\eps)$ the stopping rule determined by the discrepancy 
principle. Then there exists a constant $C$, depending on $p$ and 
$\| y \|$ only, such that
\begin{itemize}
\item[i)]  $\| \bar{x} - x_k^\eps \| \le C (\ln k)^{-p}$;
\item[ii)] $\| \gamma f_\eps - (I-T_l) x_k^\eps \| \le C k^{-1} (\ln k)^{-p}$;
\end{itemize}
for all iteration index $k$ satisfying $1 \le k \le k(\eps,f_\eps)$.
\end{proposition}

\noindent{\bf Proof.} Using the estimates in the appendices of [3] and 
[6] and the discrepancy principle one obtains
\begin{equation} \label{eq:prop-s1}
\| \bar{x} - x_k^\eps \|  \le   c\, \| y \| \, 
\big( 1 + \T\frac{1}{\mu-2} \big)\, \big( \ln(k+1) \big)^{-p} .
\end{equation}
Since $\tilde{c} := \sup_{k\in\N} \{ (\ln(k+1) / \ln(k))^{-p} \} < \infty$, 
assertion {\it i)} follows from \eqref{eq:prop-s1} with $C = c \| y \| 
(1 + \frac{1}{\mu-2})\, \tilde{c}$. To prove {\it ii)}, we again use the 
estimates in the appendices of [3] and [6] to obtain
\begin{equation} \label{eq:prop-s2}
\| \gamma f_\eps - (I-T_l) x_k^\eps \|  \le  c\, \| y \|
\big( 1 + \T\frac{2}{\mu-2} \big)\, (k-1)^{-1} \big( \ln(k+1) \big)^{-p} .
\end{equation}
Since $c^\star := \sup_{k\in\N} \{ (k+1) / k) \} < \infty$, assertion 
{\it ii)} follows from \eqref{eq:prop-s2} with $C = c \| y \| 
(1 + \frac{1}{\mu-2})\, c^\star\, \tilde{c}$.

\hfill \Box

\begin{proposition}
Set\, $k_\eps := k(\eps,f_\eps)$. Under the assumptions of the previous 
Theorem we have
\begin{itemize}
\item[i)]  $k_\eps \big( \ln(k_\eps) \big)^p \ = \ O ( \eps^{-1} )$;
\item[ii)] $\| \bar{x} - x_{k_\eps}^\eps \| \le 
             O \big(  (-\ln \sqrt{\eps})^{-p} \big)$.
\end{itemize}
\end{proposition}

\noindent{\bf Proof.} Argumenting as in the proof of Proposition~%
\ref{prop:source} we obtain
\begin{equation}
(\mu - 2) \eps \ \le \ c_1 (k_\eps - 2)^{-1} (\ln\, k_\eps)^{-p} ,
\end{equation}
from what follows
\begin{equation}
\eps^{-1} \ \ge \ c_2\, (k_\eps - 2) \, (\ln\, k_\eps)^p
    \ge  c_3\, k_\eps\, (\ln\, k_\eps)^p ,
\end{equation}
proving the first assertion. To prove {\it ii)}, we first obtain from 
the iteration rule the estimate
\begin{equation} \label{eq:prop-sc1}
\| \bar{x} - x_{k_\eps}^\eps \| \le \| F(P) v_{k_\eps} \| + \eps k_\eps ,
\end{equation}
where $P = I - T_l$. Using the estimates in the appendices of [3] and 
[6] we obtain for the first term on the right hand side of 
\eqref{eq:prop-sc1} the estimate.
\begin{equation}
\|F(P) v_{k_\eps}\|  \le  O \big( (-\ln (\eps^\frac{2}{3}))^{-p} \big) .
\end{equation}
Again using the estimates in the appendices of [3] and [6] we obtain 
for the second term on the right hand side of \eqref{eq:prop-sc1} the 
estimate
\begin{equation}
k_\eps \ = \ O \big( \eps^{-1}\, (-\ln \sqrt{\eps})^{-p} \big) .
\end{equation}
Now, substituting the last two estimates in \eqref{eq:prop-sc1}, assertion 
{\it ii)} follows.

\hfill \Box

The interest in the source conditions of logarithmic-type is motivated by 
the fact that it can be interpreted in the sense of $H^s$ regularity of 
$\bar{x} - x_1$. In order to illustrate this fact, let us consider the 
special case \,$\Omega = (-\pi,\pi)$, \,$T=1$, \,$\Lbd = (-\Delta)^{\me}$. 
We define the Sobolev spaces of periodic functions:
\begin{equation}
H^s_{per}(-\pi,\pi) \ := \ \{ \vphi(t) = \T\sum\limits_{j\in\Z} 
\vphi_j\ e^{ijt} ; \ \sum\limits_{j\in\Z} (1+j^2)^s \vphi_j^2 < \infty \}
,\ s \in \R .
\end{equation}

Clearly, the operator $T$ is well defined at \,${\cal H} = 
\overline{{\rm span} \{ \sin (jt); \ j \in \N \} }^{\|\cdot\|_{L^2}}$. 
The problem data and the initial guess can be represented in the form
\begin{equation}
f(t) = \sum_j \, f_j \sin(jt) ; \ \ \ \ 
x_1(t) = \sum_j \, \Phi_j \sin(jt) .
\end{equation}
Observe also that in ${\cal H}$ the operator $T$ can be explicitely 
represented by
\begin{equation}
 (T x_1) (t) = \sum_j \Big( (1 - \gamma e^{-j^2}) \Phi_j + 
                          \gamma f_j \Big) \sin(jt) .
\end{equation}

The next step corresponds to the choice of the parameter\, $\gamma$. Notice 
that $\bar{\lbd}$ in \eqref{eq:gammaB-wahl} satisfies 
$$ \bar{\lbd} = \inf\{ \lbd;\ \lbd \in \sigma(\Lbd) \} = 1 , $$
while the parameter $\tilde{\lbd}$ in Remark~\ref{rem:Tl} is given by
$$ \tilde{\lbd} = ( \bar{\lbd}^2 - T^{-1}\ln 2 )^{\me} = (1-\ln(2))^{\me} . $$
Thus, in this particular case, the condition \eqref{eq:gamma-wahl} for 
the choice of $\gamma$ is given by
\begin{equation}
0 < \gamma < 2\exp(\tilde{\lbd}^2T) = e
\end{equation}
and the choice $\gamma = 1$ is allowed.

From the logarithmic source condition \eqref{gl:source-cond}, with 
$y(t) = \sum_j y_j \sin(jt)$ and $q = 2p$, follows
$$ \| \bar{x} - x_1 \|_q^2 \ = \
   \T\sum_{j=1}^\infty \D(1+j^2)^{2p} \,
   \ln \Big( \frac{\exp(1)}{1 - (1-e^{-j^2})} \Big)^{-2p} y_j^2
   \ \le \ \T\sum\limits_{j=1}^\infty y_j^2
   \ \le \ \infty . $$
The reciprocal holds, i.e. if\, $\bar{x} - x_1 \in H^{2p}_{per}$, then 
exists $y \in {\cal H}$ with $\bar{x} - x_1 = f(I-{\cal T}_l)\, y$. 
Thus, the logarithmic source condition can indeed be interpreted as 
$H^p$ regularity.

\bibliographystyle{plain}
%
%
%
%
\section*{\large \bf References}
\begin{description}
\item[{[1]}] Baumeister, J.; Leit\~ao, A., {\em On iterative methods for 
solving ill-posed problems modeled by partial differential equations,} 
J. Inverse Ill-Posed Probl. {\bf 9} (2001), 13 -- 30

\item[{[2]}] Browder, F.; Petryshyn, W., {\em Construction of fixed points 
of nonlinear mappings in Hilbert space,} J. Math. Anal. Appl. {\bf 20} 
(1967), 197 -- 228

\item[{[3]}] Deuflhard, P.; Engl, H.W.; Scherzer, O., {\em A convergence 
analysis of iterative methods for the solution of nonlinear ill-posed 
problems under affinely invariant conditions,} Inverse Problems, {\bf 14} 
(1998), 1081 -- 1106

\item[{[4]}] Dotson, W.G., Jr., {\em On the Mann iterative process,} Trans. 
Am. Math. Soc. {\bf 149} (1970), 65 -- 73

\item[{[5]}] Engl, H.W.; Hanke, M.; Neubauer, A. {\em Regularization of 
Inverse Problems,} Kluwer Academic Publishers, Dordrecht, 1996 (Paperback: 
2000)

\item[{[6]}] Engl, H.W.; Leit\~ao, A., {\em A Mann iterative regularization 
method for elliptic Cauchy problems,} Numer. Funct. Anal. Optim. {\bf 22} 
(2001), 861 -- 884

\item[{[7]}] Groetsch, C., {\em A note on segmenting Mann iterates,} J. 
Math. Anal. Appl. {\bf 40} (1972), 369 -- 372

\item[{[8]}] Kirk, W.A., {\em On successive approximations for nonexpansive 
mappings in Banach spaces,} Glasgow Math. J. {\bf 12} (1971), 6 -- 9

\item[{[9]}] Mann, W., {\em Mean value methods in iteration,} Proc. Amer. 
Math. Soc. {\bf 4} (1953), 506 -- 510

\item[{[10]}] Opial, Z., {\em Weak convergence of the sequence of 
successive approximations for nonexpansive mappings,} Bull. Amer. Math. 
Soc. {\bf 73} (1967), 591 -- 597
\end{description}

\vspace{.4in}
%
%
\begin{flushleft}
Department of Mathematics \\
Federal University of Santa Catarina \\
P.O. Box 476 \\
88010-970, Florian\'opolis, SC \\
Brazil \\
{\tt aleitao@mtm.ufsc.br}   
\end{flushleft}

\end{document}